# The GAUSS LINKING INTEGRAL

# on the 3-SPHERE

# and in HYPERBOLIC 3-SPACE

*Dennis DeTurck and Herman Gluck*

We introduce here explicit integral formulas for linking, twisting, writhing and helicity on the 3-sphere and in hyperbolic 3-space.

• These formulas, like their prototypes in Euclidean 3-space, are geometric rather than just topological, in the sense that their integrands are invariant under orientation-preserving isometries of the ambient space.

• They are obtained by developing and then applying a steady-state version of classical electrodynamics in these two spaces, including an explicit Biot-Savart formula for the magnetic field and a corresponding Ampere's law contained in Maxwell's equations.

• The Biot-Savart formula leads, in turn, to upper bounds for the helicity of vector fields and lower bounds for the first eigenvalue of the curl operator on subdomains of the 3-sphere and hyperbolic 3-space.

• We give only a hint of the proofs.

June 2004

# 1. Linking integrals in $S^3$ and $H^3$.

Let $K_1 = \{x(s)\}$ and $K_2 = \{y(t)\}$ be disjoint oriented smooth closed curves in either the unit 3-sphere $S^3$ or hyperbolic 3-space $H^3$, and let $\alpha(x, y)$ denote the distance from $x$ to $y$.

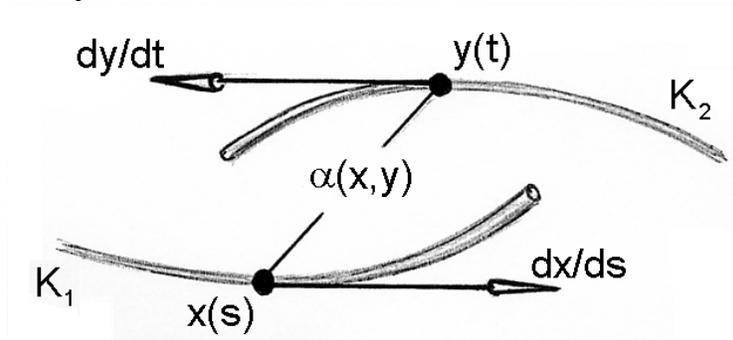

The integral formulas below contain vectors lying in different tangent spaces; in non-Euclidean settings we must move them to a common location to be combined. In $S^3$ regarded as the group of unit quaternions, the differential $L_{yx^{-1}}$ of left translation by $yx^{-1}$ moves tangent vectors from $x$ to $y$. In either $S^3$ or $H^3$, parallel transport $P_{yx}$ along the geodesic segment from $x$ to $y$ also does this.

**THEOREM 1. LINKING INTEGRALS in $S^3$ and $H^3$.**

**(1)** *On $S^3$ in left-translation format:*

$$Lk(K_1, K_2) = 1/(4\pi^2) \int_{K_1 \times K_2} L_{yx^{-1}} \, dx/ds \times dy/dt \bullet \nabla_y \, \varphi(x, y) \; ds \, dt$$
$$- \; 1/(4\pi^2) \int_{K_1 \times K_2} L_{yx^{-1}} \, dx/ds \bullet dy/dt \; ds \, dt \; ,$$

*where* $\varphi(\alpha) = (\pi - \alpha) \cot \alpha$.

**(2)** *On $S^3$ in parallel transport format:*

$$Lk(K_1, K_2) = 1/(4\pi^2) \int_{K_1 \times K_2} P_{yx} \, dx/ds \times dy/dt \bullet \nabla_y \, \varphi(x, y) \; ds \, dt \; ,$$

*where* $\varphi(\alpha) = (\pi - \alpha) \csc \alpha$.

**(3)** *On $H^3$ in parallel transport format:*

$$Lk(K_1, K_2) = 1/(4\pi) \int_{K_1 \times K_2} P_{yx} \, dx/ds \times dy/dt \bullet \nabla_y \, \varphi(x, y) \; ds \, dt \; ,$$

*where* $\varphi(\alpha) = \operatorname{csch} \alpha$.

For each of the kernel functions $\varphi(\alpha)$ above, we have written $\varphi(x, y)$ in place of $\varphi(\alpha(x, y))$.

- The kernel functions appearing in these integrals have the following significance.

    In formula (1), the function $\varphi(\alpha) = (-1/(4\pi^2))\,(\pi - \alpha)\cot\alpha$, where $\alpha$ is distance from a fixed point, is the fundamental solution of the Laplacian on $S^3$,

    $$\Delta\varphi = \delta - 1/(2\pi^2)\,.$$

    Here $\delta$ is the Dirac $\delta$-function, and, since the volume of $S^3$ is $2\pi^2$, the right-hand side has average value zero.

    In formula (2), the function $\varphi(\alpha) = (-1/(4\pi^2))\,(\pi - \alpha)\csc\alpha$ is the fundamental solution of a shifted Laplacian on $S^3$,

    $$\Delta\varphi - \varphi = \delta\,.$$

    In formula (3), the function $\varphi(\alpha) = (-1/(4\pi))\,\text{csch}\,\alpha$ is the fundamental solution of a shifted Laplacian on $H^3$,

    $$\Delta\varphi + \varphi = \delta\,.$$

- In formula (1), the first integral is the dominant term when the curves are small. At the other extreme, for orthogonal great circles, say, the first integral vanishes and the second integral takes the value $\pm 1$. In 3-space, the second integral always vanishes.

- In formula (2), it makes no difference that parallel transport from x to y on $S^3$ is ambiguous when y is the antipode $-x$ of x, because when $\alpha = \pi$ we have $\nabla_y(\varphi(\alpha(x, y))) = 0$.

- All of the integrands in formulas (1), (2) and (3) are invariant under orientation-preserving isometries of the ambient space. In formula (1), this follows from the fact that the group of left translations by unit quaternions is a **normal** subgroup of the group $SO(4)$ of all orientation-preserving isometries of $S^3$. In formulas (2) and (3), this is transparent.

- The fact that formula (1) on $S^3$ has no counterpart in $H^3$ is due to the simplicity of the group $SO(3, 1)$ of orientation-preserving isometries of $H^3$.

- Greg Kuperberg, after receiving a copy of this research announcement, wrote to us, "As it happens, I thought about this question a few years ago and I obtained one of your formulas, I think, but I never wrote it up." Subsequent correspondence shows that he did indeed obtain, by a beautiful geometric argument totally different from ours, an expression equivalent to formula (2).

## 2. Historical background.

*Carl Friedrich Gauss*, in a half-page paper dated January 22, 1833, gave an integral formula for the linking number of two disjoint smooth closed curves in Euclidean 3-space:

$$Lk(K_1, K_2) \;=\; 1/(4\pi) \int_{K_1 \times K_2} \; dx/ds \times dy/dt \bullet (x - y)/|x - y|^3 \; ds \, dt \, .$$

According to historian Moritz Epple (1998), Gauss was interested in computing the linking number of the earth's orbit with the orbits of certain asteroids, and although he presented his formula without proof, it is believed that he simply counted up how many times the vector from the earth to the asteroid covered the "heavenly sphere" ..... a degree-of-map argument. Gauss undoubtedly knew another proof: run a current through the first loop, and calculate the circulation of the resulting magnetic field around the second loop. By Ampere's Law, this circulation is equal to the total current enclosed by the second loop, which means the current flowing along the first loop, multiplied by the linking number of the two loops. Then the Biot-Savart formula (1820) for the magnetic field leads directly to Gauss's linking integral.

*Georges Calugareanu*, in a series of three papers (1959-1961), let the two curves in Gauss's linking integral come together, and defined a real-valued invariant, later called the "writhing number" by F. Brock Fuller (1971), of a smooth simple closed curve in Euclidean 3-space, which measures the extent to which the curve wraps and coils around itself:

$$Wr(K) \;=\; 1/(4\pi) \int_{K \times K} \; dx/ds \times dy/dt \bullet (x - y)/|x - y|^3 \; ds \, dt \, .$$

In these papers, Calugareanu also discovered the formula LINK = TWIST + WRITHE, in which LINK is the linking number of the two edges of a closed ribbon, TWIST measures the extent to which the ribbon twists around one of its edges, and WRITHE is the writhing number of that edge.

Calugareanu proved this formula under the assumption that the simple closed curve K has nowhere-vanishing curvature although, as we note in section 11, the basic ideas for proving the formula without this assumption are already present in his papers. Calugareanu's formula without this restriction was proved by James White (1969) in his thesis, using a totally different approach based on ideas of William Pohl (1968a, b).

***Lodewijk Woltjer*** introduced in 1958 the notion of "helicity" of a vector field V defined on a domain $\Omega$ in Euclidean 3-space,

$$H(V) = 1/(4\pi) \int_{\Omega \times \Omega} V(x) \times V(y) \bullet (x - y)/|x - y|^3 \; d(vol_x) \, d(vol_y) \, ,$$

as an invariant during ideal MHD evolution of plasma fields. Keith Moffatt (1969), recognizing that this quantity measures the extent to which the field lines of V wrap and coil around one another, named it "helicity" and showed that Woltjer's original formula could be written in the above form.

The formulas of Gauss and Calugareanu serve as cornerstones of geometric knot theory in Euclidean space, with rich connections to molecular biology, especially in the study of knotted duplex DNA and of the enzymes which affect it, and, through the related formula of Woltjer and Moffatt, to fluid dynamics, electrodynamics and plasma and solar physics.

## 3. Organization of this report.

Gauss's degree-of-map derivation of his linking integral does not work on the 3-sphere $S^3$ because the set of ordered pairs of distinct points in $S^3$ deformation retracts to a 3-sphere rather than to a 2-sphere, as it does in $R^3$. But a degree-of-map derivation, suitably modified, does work in hyperbolic 3-space.

So we begin by developing a steady-state version of classical electrodynamics on $S^3$ and $H^3$ in sections 4 through 8, including an explicit Biot-Savart formula for the magnetic field and a corresponding Ampere's law contained in Maxwell's equations. In section 9, this enables us to follow Gauss's second proof of the linking integrals.

With the linking integrals in hand, in sections 10 and 11 we follow the lead of Calugareanu towards the writhing and twisting integrals and the formula Link = Twist + Writhe.

In section 12 we follow the lead of Woltjer towards the helicity integrals.

In section 13 we derive upper bounds for the helicity of vector fields in $R^3$, $S^3$ and $H^3$, and in section 14 we derive lower bounds for the first eigenvalue of the curl operator on subdomains of these spaces.

For a development of electrodynamics on bounded subdomains of the 3-sphere, and for applications of the linking, writhing and helicity integrals in this setting, including a proposed application to solar physics, see the Ph.D. thesis of Jason Parsley (2004).

We take pleasure in acknowledging Jason Parsley's help throughout the preparation of this report, and especially in the development of electrodynamics on $S^3$.

# 4. Magnetic fields on $S^3$ and $H^3$.

In Euclidean 3-space $R^3$, the classical convolution formula of Biot and Savart gives the magnetic field BS(V) generated by a compactly supported current flow V:

$$BS(V)(y) = 1/(4\pi) \int_{R^3} V(x) \times (y - x) / |y - x|^3 \, dx.$$

For simplicity, we write dx to mean $d(vol_x)$.

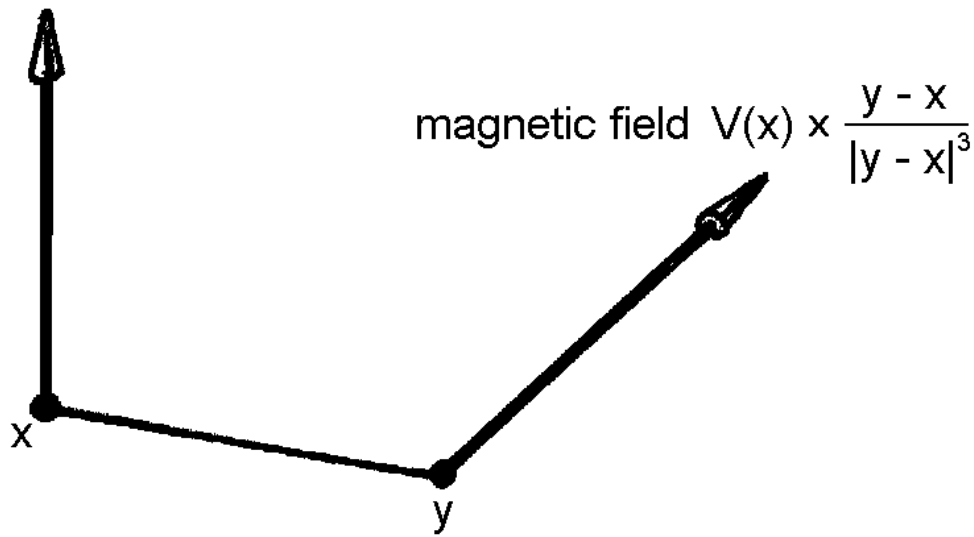

The function $\varphi(r) = -1/(4\pi r)$ is the fundamental solution of the Laplacian in $R^3$,

$$\Delta\varphi = \delta.$$

Because the gradient of $\varphi$ at y,

$$\nabla_y \varphi(x, y) = (1/(4\pi)) (y - x)/|y - x|^3,$$

is the kernel function appearing in the Biot-Savart formula above, this formula can be written as

$$BS(V)(y) = \int_{R^3} V(x) \times \nabla_y \varphi(x, y) \, dx.$$

In $R^3$, the magnetic field BS(V) of a smooth, compactly supported current flow V has four important properties:

**(1) It has finite energy, that is $\int_{R^3}$ BS(V)(y) • BS(V)(y) dy $< \infty$.**

**(2) It is divergence-free, $\nabla$ • BS(V) = 0.**

**(3) If V is a gradient field, then BS(V) = 0.**

**(4) If V is divergence-free, then $\nabla \times$ BS(V) = V.**

Properties (1), (2) and (4) are well known, and the simple verification of (3) can be found in our paper (2001).

It is easy to see that in $R^3$ there can't be two different linear operators BS with the above properties, so in this familiar setting we have both existence and uniqueness.

We take these properties as constraints on the definition of the Biot-Savart operator in $S^3$ and $H^3$, giving us uniqueness there as well.

The following theorem asserts existence, and provides explicit integral formulas.

# THEOREM 2. BIOT-SAVART INTEGRALS in $S^3$ and $H^3$.

*The Biot-Savart operator exists in $S^3$ and $H^3$, and is given by the following formulas.*

**(1)** *On $S^3$ in left-translation format:*

$$BS(V)(y) = \int_{S^3} L_{yx^{-1}} V(x) \times \nabla_y \varphi_0(x, y) \, dx$$
$$- 1/(4\pi^2) \int_{S^3} L_{yx^{-1}} V(x) \, dx$$
$$+ 2 \nabla_y \int_{S^3} L_{yx^{-1}} V(x) \bullet \nabla_y \varphi_1(x, y) \, dx,$$

*where* $\varphi_0(\alpha) = (-1/(4\pi^2))(\pi-\alpha) \cot \alpha$ *and* $\varphi_1(\alpha) = (-1/(16\pi^2)) \alpha (2\pi-\alpha)$.

**(2)** *On $S^3$ in parallel transport format:*

$$BS(V)(y) = \int_{S^3} P_{yx} V(x) \times \nabla_y \varphi(x, y) \, dx,$$

*where* $\varphi(\alpha) = (-1/(4\pi^2))(\pi - \alpha) \csc \alpha$.

**(3)** *On $H^3$ in parallel transport format:*

$$BS(V)(y) = \int_{H^3} P_{yx} V(x) \times \nabla_y \varphi(x, y) \, dx,$$

*where* $\varphi(\alpha) = (-1/(4\pi)) \operatorname{csch} \alpha$.

• In formula (1), the function $\varphi_1(\alpha) = (-1/(16\pi^2)) \alpha (2\pi-\alpha)$ satisfies the equation

$$\Delta \varphi_1 = \varphi_0 - [\varphi_0],$$

where $[\varphi_0]$ denotes the average value of $\varphi_0$ over $S^3$. The other kernel functions already appeared in the linking integrals in Theorem 1.

- A vector field V on $S^3$ is said to be *left-invariant* if

$$L_{yx^{-1}} V(x) = V(y) \qquad \text{for all } x, y \in S^3.$$

For example, using multiplication of unit quaternions, the vector fields

$$V_1(x) = x\,i, \quad V_2(x) = x\,j, \quad V_3(x) = x\,k$$

form a pointwise orthonormal basis for the 3-dimensional subspace of left-invariant vector fields on $S^3$.

Left-invariant and right-invariant vector fields on $S^3$ are tangent to the orbits of Hopf fibrations of $S^3$ by parallel great circles.

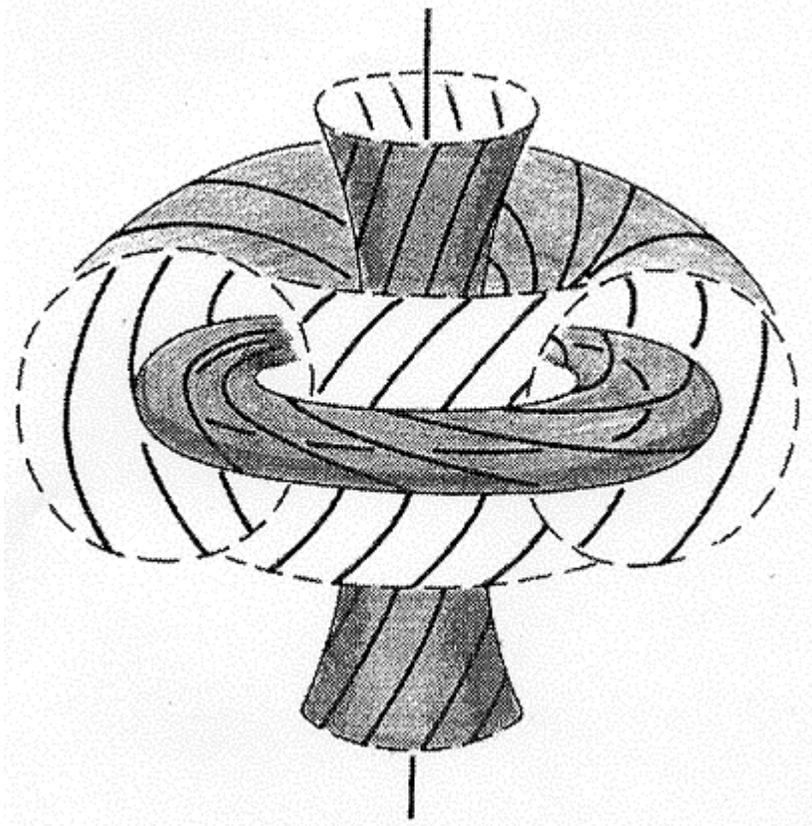

It is not surprising that left-invariant fields play an important role in the left-translation versions of our theorems. In particular, if V is a left-invariant field on $S^3$, then V is divergence-free and is a curl eigenfield with eigenvalue $-2$,

$$\nabla \times V = -2\,V\,.$$

Moreover, the Biot-Savart operator is the inverse of curl on divergence-free fields on $S^3$,

$$BS(V) = -\tfrac{1}{2}\,V\,.$$

In formula (1) of Theorem 2,

- the first integral is nonzero unless V is a left-invariant field;

- the second integral is nonzero if and only if V has a nonzero component, with respect to the $L^2$ inner product $<V, W> = \int_{S^3} V \bullet W\ d(vol)$, in the direction of the 3-dimensional subspace of left-invariant fields;

- the third integral is zero if and only if V is divergence-free;

- if V is a gradient field, then the first and third integrals are nonzero and cancel each other, while the second integral is zero.

## 5. Proof scheme for Theorem 2, formula (1).

Let $VF(S^3)$ denote the set of all smooth ($C^\infty$) vector fields on $S^3$, regarded as an infinite-dimensional vector space with the $L^2$ inner product

$$<V, W> = \int_{S^3} V \bullet W \, d(vol) .$$

The *vector Laplacian*,

$$\Delta V = - \nabla \times \nabla \times V + \nabla(\nabla \bullet V) ,$$

corresponds to the Laplacian on forms,

$$\Delta \beta_V = - *d*d \, \beta_V + d*d* \, \beta_V ,$$

where the one-form $\beta_V$ is defined by $\beta_V(W) = V \bullet W$.

The vector Laplacian $\Delta: VF(S^3) \to VF(S^3)$ is an unbounded self-adjoint linear operator which is bijective, and which takes divergence-free fields to divergence-free fields, and gradient fields to gradient fields.

The *Green's operator* $Gr: VF(S^3) \to VF(S^3)$ is defined to be the inverse of the vector Laplacian. It is a bounded self-adjoint linear operator which takes divergence-free fields to divergence-free fields and gradient fields to gradient fields.

Let $V$ be a smooth vector field on $S^3$. Thinking of $V$ as a steady current flow, we define the corresponding magnetic field, and the Biot-Savart operator $BS$, by

$$BS(V) = - \nabla \times Gr(V) .$$

It is easy to check that the magnetic field $BS(V)$ satisfies the four properties required of it in the preceding section.

Once we get an explicit formula for the Green's operator $Gr(V)$, the negative of its curl will lead to the desired explicit formula for $BS(V)$. This in turn will lead us in section 9 to the linking integral formula given in Theorem 1, part (1).

The fact that an explicit formula for the Green's operator can lead (via a different route) to an explicit formula for the linking integral was known to de Rham (1960, section 33). The issue is to find the formula for the Green's operator, and we turn now to showing how to do this on $S^3$.

Working in left-translation format, it is natural to hope, by analogy with $R^3$, that

$$A(V, \varphi)(y) = \int_{S^3} L_{yx^{-1}} V(x) \varphi(x, y) \, dx$$

will be the Green's operator on $S^3$ when we choose $\varphi$ to be the fundamental solution of the Laplacian. But it doesn't work out this way.

However, $\varphi$ is not the only possible kernel function, and $A(V, \varphi)$ is not the only way of convolving a vector field with a kernel to get another vector field.

Consider the convolutions

$$A(V, \varphi)(y) = \int_{S^3} L_{yx^{-1}} V(x) \varphi(x, y) \, dx$$

$$B(V, \varphi)(y) = \int_{S^3} L_{yx^{-1}} V(x) \times \nabla_y \varphi(x, y) \, dx$$

$$G(V, \varphi)(y) = \nabla_y \int_{S^3} L_{yx^{-1}} V(x) \bullet \nabla_y \varphi(x, y) \, dx \, .$$

In $R^3$, if we let $\varphi$ be the fundamental solution of the scalar Laplacian, then the convolution $A(V, \varphi)$ gives the Green's operator, $B(V, \varphi)$ gives the magnetic field, and $G(V, \varphi)$ gives the $\partial E/\partial t$ term in Maxwell's equation.

On $S^3$, it turns out that a linear combination of all three types of vector convolution, each with a different kernel function $\varphi$, will give the Green's operator.

To prove this, we develop a "calculus of vector convolutions" in left-translation format, in which we compute formulas for curls, divergences and Laplacians of convolutions of types $A$, $B$ and $G$. The key formulas are the following:

$$\Delta A(V, \varphi) = A(V, \Delta\varphi) \quad - 4 A(V, \varphi) \quad - 2 B(V, \varphi)$$

$$\Delta B(V, \varphi) = B(V, \Delta\varphi) \quad + 2 A(V, \Delta\varphi) - 2 G(V, \varphi)$$

$$\Delta G(V, \varphi) = G(V, \Delta\varphi) \, .$$

The large spaces after the first terms on the right hand sides of the $A$ and $B$ lines serve as a reminder that in Euclidean space $R^3$, only those first terms appear. In particular, in $R^3$ we have $\Delta A(V, \varphi) = A(V, \Delta\varphi)$. So naturally, when we choose $\varphi$ to be a fundamental solution of the scalar Laplacian in $R^3$, we get $\Delta A(V, \varphi) = A(V, \delta) = V$, which tells us that $A(V, \varphi) = Gr(V)$ there.

On the 3-sphere $S^3$, we get the Green's operator by adding one copy of line A to two copies of line B and four copies of line G, each with a different choice of kernel function $\varphi$, as follows.

**THEOREM 3.** *On $S^3$ in left-translation format, the Green's operator on vector fields is given by*

$$Gr(V)(y) = \int_{S^3} L_{yx^{-1}} V(x) \, \varphi_0(x, y) \, dx$$

$$+ \, 2 \int_{S^3} L_{yx^{-1}} V(x) \times \nabla_y \varphi_1(x, y) \, dx$$

$$+ \, 4 \nabla_y \int_{S^3} L_{yx^{-1}} V(x) \bullet \nabla_y \varphi_2(x, y) \, dx \, ,$$

*where*

$$\varphi_0(\alpha) = (-1/(4\pi^2)) \, (\pi - \alpha) \cot \alpha$$

$$\varphi_1(\alpha) = (-1/(16\pi^2)) \, \alpha \, (2\pi - \alpha)$$

$$\varphi_2(\alpha) = (-1/(192\pi^2)) \, \big(3\alpha(2\pi - \alpha) + 2\alpha(\pi - \alpha)(2\pi - \alpha) \cot \alpha\big) \, ,$$

*and these three kernel functions are related by*

$$\varphi_2 \xrightarrow{\Delta} \varphi_1 - [\varphi_1] \xrightarrow{\Delta} \varphi_0 - [\varphi_0] \xrightarrow{\Delta} \delta - [\delta] \, .$$

The negative curl of the above formula for the Green's operator,

$$BS(V) = -\nabla \times Gr(V) \, ,$$

leads to formula (1) of Theorem 2 for the Biot-Savart operator on $S^3$.

# 6. Proof scheme for Theorem 2, formula (2).

Now we are working on $S^3$ in parallel transport format, and use an entirely different approach.

**KEY LEMMA, spherical version.**

$$\nabla_y \times \{P_{yx}V(x) \times \nabla_y \varphi\} - \nabla_y \{V(x) \bullet \nabla_x (\cos \alpha \ \varphi)\}$$
$$= (\Delta \varphi - \varphi)(V(x) - <V(x), y> y).$$

The above lemma is quite flexible: $V$ can be a smooth vector field defined on all of $S^3$, or just on a bounded subdomain $\Omega$ of $S^3$, or just along a smooth closed curve $K$ in $S^3$, or even just at the single point $x$. The kernel function $\varphi(\alpha)$ is any smooth function, which typically blows up at $\alpha = 0$ and is asymptotic there to $-1/(4\pi\alpha)$. As usual, $\alpha = \alpha(x, y)$ is the distance on $S^3$ between $x$ and $y$. The dot product $\bullet$ on the left hand side of the above equation is in the tangent space to $S^3$ at $x$, while the inner product $<,>$ on the right hand side is in $R^4$.

The proof of the lemma exploits the interplay between the geometry of $S^3$ and the linear algebra of the space $R^4$ in which it lies. Here are four samples:

- We define a triple vector product $[A, B, C]$ in $R^4$, whose value is a vector orthogonal to the 3-space spanned by $A$, $B$ and $C$, and whose length is the 3-dimensional volume of the parallelepiped spanned by these three vectors.

- If $u$ and $v$ are both tangent to $S^3$ at $x$, then $u \times v = [x, u, v]$.

- If $x$ and $y$ are points of $S^3$ and $V(x)$ is tangent to $S^3$ at $x$, then

$$\nabla_y \times [x, V(x), y] = 2 <x, y> V(x) - 2 <V(x), y> x.$$

- Under these same hypotheses,

$$[y, x, [y, x, V(x)]] = -\sin^2\alpha \ V(x) + <V(x), y> y - \cos \alpha <V(x), y> x.$$

Once the Key Lemma is proved, we apply it using the kernel function

$$\varphi(\alpha) = (-1/(4\pi^2))(\pi - \alpha)\csc\alpha,$$

which, as we noted earlier, is a fundamental solution of the shifted Laplacian on $S^3$,

$$\Delta\varphi - \varphi = \delta.$$

With this choice of kernel function, we have

$$(\cos\alpha)\varphi(\alpha) = \varphi_0(\alpha) = (-1/(4\pi^2))(\pi - \alpha)\cot\alpha,$$

which is the fundamental solution of the Laplacian on $S^3$,

$$\Delta\varphi_0 = \delta - 1/(2\pi^2).$$

With these choices, the Key Lemma reads

$$\nabla_y \times \{P_{yx}V(x) \times \nabla_y\varphi\} - \nabla_y\{V(x) \bullet \nabla_x \varphi_0\}$$

$$= \delta(x, y)(V(x) - <V(x), y> y).$$

Keeping these choices of kernel functions, we now define the Biot-Savart operator by formula (2) of Theorem 2,

$$BS(V)(y) = \int_{S^3} P_{yx}V(x) \times \nabla_y\varphi(x, y)\, dx.$$

Then we integrate the Key Lemma over $S^3$ and get

$$\int_{S^3} \nabla_y \times \{P_{yx}V(x) \times \nabla_y\varphi\}\, dx - \int_{S^3} \nabla_y\{V(x) \bullet \nabla_x \varphi_0\}\, dx$$

$$= \int_{S^3} \delta(x, y)(V(x) - <V(x), y> y)\, dx.$$

The integral on the right hand side above equals

$$V(y) - <V(y), y> y = V(y),$$

since $V(y)$ is a tangent vector to $S^3$ at $y$, and hence orthogonal to $y$ in $R^4$.

Taking the curl and the gradient outside the two integrals on the left-hand side, we then get

$$\nabla_y \times \int_{S^3} P_{yx} V(x) \times \nabla_y \varphi(x, y) \, dx - \nabla_y \int_{S^3} V(x) \bullet \nabla_x \varphi_0(x, y) \, dx = V(y),$$

or

$$\nabla_y \times BS(V)(y) = V(y) + \nabla_y \int_{S^3} V(x) \bullet \nabla_x \varphi_0(x, y) \, dx.$$

Integrating by parts over $S^3$, we have

$$\int_{S^3} V(x) \bullet \nabla_x \varphi_0(x, y) \, dx = -\int_{S^3} \nabla_x \bullet V(x) \, \varphi_0(x, y) \, dx,$$

and inserting this, we get

$$\nabla_y \times BS(V)(y) = V(y) - \int_{S^3} \nabla_x \bullet V(x) \, \varphi_0(x, y) \, dx.$$

If we think of the vector field $V(x)$ as a steady current, then the negative divergence, $-\nabla_x \bullet V(x)$, is the time rate of accumulation of charge at $x$, and hence the integral

$$-\int_{S^3} \nabla_x \bullet V(x) \, \varphi_0(x, y) \, dx$$

is the time rate of increase of the electric field at $y$ (see section 8). Thus the above equation is Maxwell's equation

$$\nabla \times B = J + \partial E / \partial t.$$

With this line of reasoning in mind, we view the Key Lemma, with our specific choice of kernel functions, as a pointwise version of Maxwell's equation.

It is now easy to verify that $BS(V)$ satisfies the four properties required of it in section 4.

**Remark.** Jason Parsley, in his thesis (2004), proves Maxwell's equation in the case that the current flow $J$ runs only within some bounded subdomain $\Omega \subset S^3$. In that situation, the $\partial E/\partial t$ term also reflects the time rate of change of the electric field due to the gradual deposit or removal of charge along the boundary of $\Omega$. His proof of this in left-invariant format inspired our own arguments.

## 7. Proof scheme for Theorem 2, formula (3).

We replace $R^4$ by Minkowski space $R^{1,3}$ with the indefinite inner product

$$<x, y> \; = \; x_0 y_0 \; - \; x_1 y_1 \; - \; x_2 y_2 \; - \; x_3 y_3 \; ,$$

and then regard hyperbolic 3-space as the hyperboloid

$$H^3 \; = \; \{ \, x \in R^{1,3} : \; <x, x> \; = \; 1 \text{ and } x_0 > 0 \, \} \; \subset \; R^{1,3} \; .$$

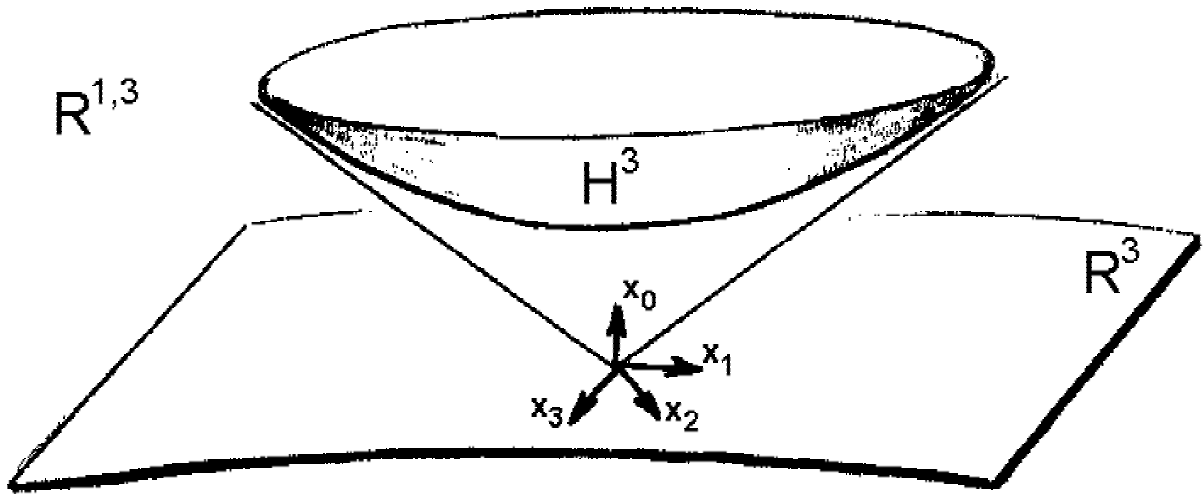

If $u$ and $v$ are tangent to $H^3$ at $x$, we define the induced inner product $\bullet$ on $T_x H^3$ by

$$u \bullet v \; = \; - <u, v> \; .$$

This time we exploit the interplay between the geometry of $H^3$ and the linear algebra of $R^{1,3}$ to prove the following

**KEY LEMMA, hyperbolic version.**

$$\nabla_y \times \{P_{yx} V(x) \times \nabla_y \varphi\} \; - \; \nabla_y \{V(x) \bullet \nabla_x (\cosh \alpha \; \varphi)\}$$
$$= \; (\Delta \varphi \; + \; \varphi) \, (V(x) \; - \; <V(x), y> y) \; .$$

Comparing the Key Lemmas in $S^3$ and $H^3$, we see that the only changes are that $\cos \alpha$ has been replaced by $\cosh \alpha$ on the left-hand side, and $\Delta \varphi - \varphi$ has been replaced by $\Delta \varphi + \varphi$ on the right-hand side.

After proving the Key Lemma for $H^3$, we apply it using the kernel function

$$\varphi(\alpha) = (-1/(4\pi)) \operatorname{csch} \alpha,$$

which, as we noted earlier, is a fundamental solution of the shifted Laplacian on $H^3$,

$$\Delta \varphi + \varphi = \delta.$$

With this choice of kernel function, we have

$$(\cosh \alpha)\, \varphi(\alpha) = \varphi_0(\alpha) = (-1/(4\pi)) \coth \alpha,$$

which is the fundamental solution of the Laplacian on $H^3$,

$$\Delta \varphi_0 = \delta.$$

With these choices, the Key Lemma reads

$$\nabla_y \times \{P_{yx} V(x) \times \nabla_y \varphi\} - \nabla_y \{V(x) \bullet \nabla_x \varphi_0\}$$
$$= \delta(x, y)\,(V(x) - <V(x), y> y).$$

When we define the Biot-Savart operator by formula (3) of Theorem 2,

$$BS(V)(y) = \int_{H^3} P_{yx} V(x) \times \nabla_y \varphi(x, y)\, dx,$$

and follow the line of reasoning from the previous section, we obtain Maxwell's equation

$$\nabla \times B = J + \partial E/\partial t,$$

and learn that $BS(V)$ satisfies the four properties required of it in section 4.

# 8. Maxwell's equations on $S^3$ and $H^3$.

Maxwell's equations, so familiar to us in Euclidean 3-space $R^3$, can also be shown to hold on $S^3$ and $H^3$, with the appropriate definitions, as follows.

Let $\rho$ be a smooth real-valued function, which we think of as a charge density. On $R^3$ and $H^3$, we ask that $\rho$ have compact support, while on $S^3$ we ask that it have average value zero.

In each case, we define the corresponding electric field $E = E(\rho)$ by the formula

$$E(\rho)(y) = \nabla_y \int \rho(x) \varphi(x, y) \, dx ,$$

where $\varphi$ is the fundamental solution of the Laplacian.

Let J be a smooth vector field, which we think of as a steady-state (i.e., time-independent) current distribution. On $R^3$ and $H^3$, we ask that J have compact support, while on $S^3$ we impose no restriction on J.

In each case, we define the corresponding magnetic field $B = BS(J)$ by one of the formulas in Theorem 2. Since J is steady-state, so is B, and therefore $\partial B/\partial t = 0$.

## THEOREM 4. MAXWELL'S EQUATIONS in $S^3$ and $H^3$.

*With these definitions, Maxwell's equations,*

$$\nabla \cdot E = \rho \qquad \nabla \times E = 0$$
$$\nabla \cdot B = 0 \qquad \nabla \times B = J + \partial E/\partial t ,$$

*hold on $S^3$ and in $H^3$, just as they do in $R^3$.*

The proofs of Maxwell's equations are either self-evident or immediate consequences of information from previous sections, as follows.

We have

$$\nabla_y \bullet E(\rho)(y) = \nabla_y \bullet \nabla_y \int \rho(x) \varphi(x, y) \, dx$$

$$= \Delta_y \int \rho(x) \varphi(x, y) \, dx$$

$$= \rho(y) \, ,$$

because $\varphi$ is the fundamental solution of the Laplacian.

We have $\nabla \times E = 0$ because $E$ is a gradient field.

We have $\nabla \bullet B = 0$ because that is one of the four properties required of $B$.

We have already noted in sections 6 and 7 how we get the last Maxwell equation, $\nabla \times B = J + \partial E/\partial t$, from the Key Lemma.

• The equation $d\rho/dt = - \nabla \bullet J$ expresses the fact that the divergence of the current flow $J$ corresponds to the time rate of accumulation or dissipation of charge, which then leads to the $\partial E/\partial t$ term in Maxwell's equation above, as we noted in section 6.

• If the current flow $J$ is divergence-free, then $\partial E/\partial t = 0$, and hence $\nabla \times B = J$. In that case, we get Ampere's Law (just the curl theorem), which says that the circulation of $B$ around a loop equals the flux of $J$ through any surface bounded by that loop.

## 9. Proof scheme for Theorem 1.

With steady-state electrodynamics in hand, we can follow what we imagine to be Gauss's second line of reasoning and obtain the linking integral formulas on $S^3$ and $H^3$ given in Theorem 1.

To do this, run a unit current $V = dx/ds$ along the first loop $K_1$, and use one of the three Biot-Savart formulas from Theorem 2 to calculate the resulting magnetic field $B$. For example, on $S^3$ in parallel transport format, we have

$$B(y) = BS(dx/ds)(y) = \int_{K_1} P_{yx} \, dx/ds \times \nabla_y \varphi(x, y) \, ds \, ,$$

where $\varphi(\alpha) = (-1/(4\pi^2)) \, (\pi - \alpha) \csc \alpha$, and where the integration takes place with respect to arc length along $K_1$ rather than with respect to volume in $S^3$.

We intend to calculate the circulation of this magnetic field $B$ along the second loop $K_2$, and then apply Ampere's Law to get the formula for the linking integral.

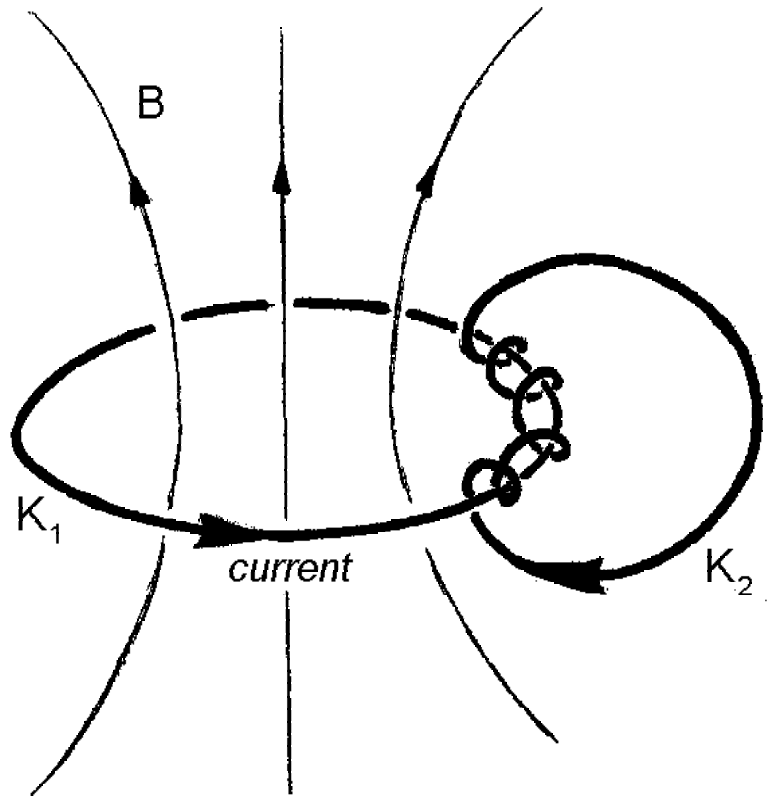

We can view $V$ as the limit of divergence-free current flows $V_\varepsilon$ supported in smaller and smaller tubular neighborhoods $N_\varepsilon$ of $K_1$, each disjoint from $K_2$, and each with unit flux through a cross-sectional disk. Let $B_\varepsilon$ be the corresponding magnetic field, defined throughout $S^3$ or $H^3$. Since $V_\varepsilon$ is divergence-free, we have $\nabla \times B_\varepsilon = V_\varepsilon$, and so by Ampere's Law, the circulation of $B_\varepsilon$ around the second loop $K_2$ equals the unit flux of $V_\varepsilon$ through the cross-sectional disks of $N_\varepsilon$, multiplied by the linking number of $K_1$ and $K_2$.

Now the magnetic fields $B_\varepsilon$ converge uniformly to $B$ on the complement of any fixed tubular neighborhood of $K_1$, and in particular along the loop $K_2$. Hence the circulation of $B$ along $K_2$ also equals the linking number of $K_1$ and $K_2$.

The three formulas for the magnetic field given in Theorem 2 then lead quickly to the three linking integrals in Theorem 1.

For example, working on $S^3$ in parallel transport format,

$$Lk(K_1, K_2) = \text{circulation of } B \text{ along } K_2$$

$$= \int_{K_2} BS(dx/ds)(y) \bullet dy/dt \, dt$$

$$= \int_{K_2} \left( \int_{K_1} P_{yx} \, dx/ds \times \nabla_y \varphi(x, y) \, ds \right) \bullet dy/dt \, dt$$

$$= \int_{K_1 \times K_2} P_{yx} \, dx/ds \times \nabla_y \varphi(x, y) \bullet dy/dt \, ds \, dt$$

$$= -\int_{K_1 \times K_2} P_{yx} \, dx/ds \times dy/dt \bullet \nabla_y \varphi(x, y) \, ds \, dt \,,$$

where $\varphi(\alpha) = (-1/(4\pi^2))(\pi - \alpha) \csc \alpha$,

$$= 1/(4\pi^2) \int_{K_1 \times K_2} P_{yx} \, dx/ds \times dy/dt \bullet \nabla_y \varphi(x, y) \, ds \, dt \,,$$

where now $\varphi(\alpha) = (\pi - \alpha) \csc \alpha$, giving us formula (2) of Theorem 1.

• Note that the third integral in the formula for $BS(V)$ in left-translation format on $S^3$, given in Theorem 2 part (1), is missing from the corresponding linking integral, given in Theorem 1 part (1), because the approximating current flows $V_\varepsilon$ are all divergence-free.

# 10. Link, twist and writhe on $S^3$ and $H^3$.

We follow Calugareanu's lead and replace the two closed curves $K_1$ and $K_2$ in our linking integrals by a single closed curve $K$, which we must assume to be *simple*. Now the points $x(s)$ and $y(t)$ run along the same curve, and therefore can coincide, making three out of the four linking integrals appearing in Theorem 1 seem improper (the second integral on $S^3$ in left-translation format is clearly proper). But as Calugareanu noted for the corresponding $R^3$ integral, the three threatening integrands actually vanish along the diagonal of $K \times K$, so the integrals converge. The resulting value $Wr(K)$ is a real number called the *writhing number* of $K$. It is, as we mentioned earlier, a measure of the extent to which $K$ wraps and coils around itself.

**DEFINITION of the WRITHING INTEGRALS in $S^3$ and $H^3$.**

**(1)** *On $S^3$ in left-translation format:*

$$Wr_L(K) =_{defn} 1/(4\pi^2) \int_{K \times K} L_{yx^{-1}} \, dx/ds \times dy/dt \bullet \nabla_y \varphi(x, y) \, ds \, dt$$
$$- 1/(4\pi^2) \int_{K \times K} L_{yx^{-1}} \, dx/ds \bullet dy/dt \, ds \, dt \, ,$$

*where* $\varphi(\alpha) = (\pi - \alpha) \cot \alpha$ .

**(2)** *On $S^3$ in parallel transport format:*

$$Wr_P(K) =_{defn} 1/(4\pi^2) \int_{K \times K} P_{yx} \, dx/ds \times dy/dt \bullet \nabla_y \varphi(x,y) \, ds \, dt \, ,$$

*where* $\varphi(\alpha) = (\pi - \alpha) \csc \alpha$ .

**(3)** *On $H^3$ in parallel transport format:*

$$Wr_P(K) =_{defn} 1/(4\pi) \int_{K \times K} P_{yx} \, dx/ds \times dy/dt \bullet \nabla_y \varphi(x,y) \, ds \, dt \, ,$$

*where* $\varphi(\alpha) = \operatorname{csch} \alpha$ .

The two versions of the writhing number on $S^3$ are *not* the same, and one can show that
$$Wr_L(K) = Wr_P(K) + L/(2\pi) \, ,$$

where $L$ is the length of the closed curve $K$. The parallel transport version of writhe is the more intuitively satisfying, since in this version the writhing number of a great circle is zero.

We turn next to the definition of "twist".

Let $K$ be a smooth simple closed curve in $S^3$ or $H^3$, parametrized by arc length $s$, with $0 \leq s \leq L$. Let $x(s)$ be a moving point along $K$, and $T(s) = x'(s)$ the unit tangent vector field.

Let $V(s)$ be a unit normal vector field along $K$. Our intention is to define the (total) twist of $V$ along $K$ by a formula such as

$$Tw(V) = 1/(2\pi) \int_K T(s) \times V(s) \bullet V'(s) \, ds,$$

as in Euclidean 3-space.

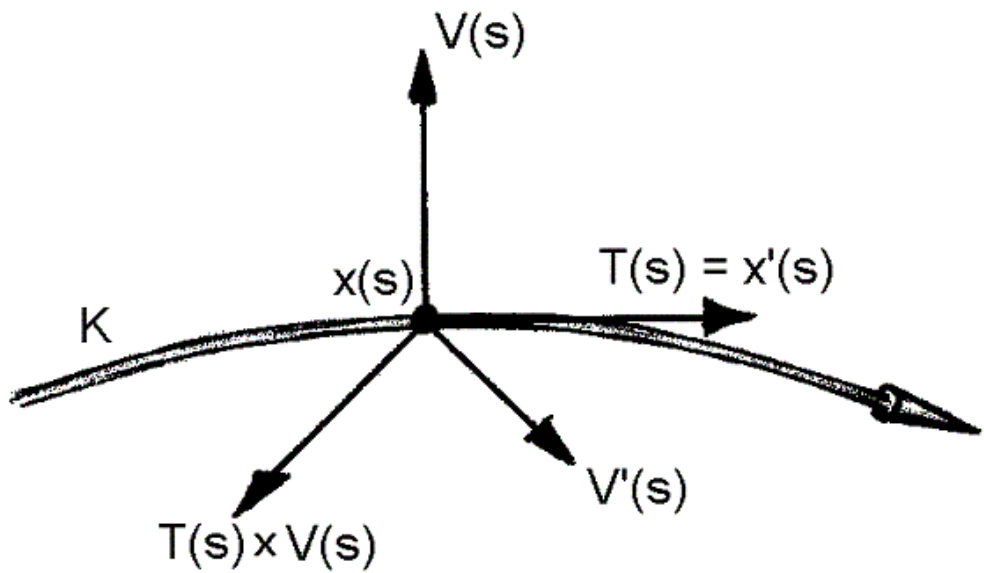

But on $S^3$ there are two flavors of twist, according as $V'(s)$ is calculated as a "left-invariant" derivative or a covariant derivative. If we fall back into Euclidean mode and write

$$V'(s) = \lim_{h \to 0} (V(s+h) - V(s))/h,$$

then the vectors $V(s+h)$ and $V(s)$ lie in different tangent spaces, and we must move them together in order to subtract. If we use left-translation in the group $S^3$ to move $V(s+h)$ back to the tangent space at $x(s)$ which contains $V(s)$, then the resulting limit is the *left-invariant derivative* $V'_L(s)$. If we use parallel transport to move $V(s+h)$ back, then the resulting limit is the *covariant derivative* $V'_P(s)$.

The two flavors of twist on $S^3$ are then given by

(1) $$Tw_L(V) \;=\; 1/(2\pi) \int_K \; T(s) \times V(s) \bullet V'_L(s) \; ds$$

and

(2) $$Tw_P(V) \;=\; 1/(2\pi) \int_K \; T(s) \times V(s) \bullet V'_P(s) \; ds \;.$$

It is an easy exercise to prove that

$$Tw_L(V) \;=\; Tw_P(V) \;-\; L/(2\pi) \;,$$

where L is the length of the curve K.

***Example.*** Consider the great circle $K = \{(\cos s, \sin s, 0, 0): 0 \le s \le 2\pi\}$ on $S^3$, and along it the unit normal vector field $V(s) = (0, 0, \cos s, \sin s)$. Then we have $Tw_L(V) = 0$ and $Tw_P(V) = 1$.

In hyperbolic 3-space $H^3$, we have only the parallel transport version of twist,

(3) $$Tw_P(V) \;=\; 1/(2\pi) \int_K \; T(s) \times V(s) \bullet V'_P(s) \; ds \;.$$

Now consider in $S^3$ or $H^3$ a narrow ribbon of width $\varepsilon$ obtained by starting with a simple closed curve $K = \{x(s)\}$ and then exponentiating a unit normal vector field $V(s)$ along $K$. One edge of this ribbon is the original curve $K$, and the other edge is the curve $K_\varepsilon = \{y_\varepsilon(s)\}$, given explicitly by

$$y_\varepsilon(s) = \cos\varepsilon\ x(s) + \sin\varepsilon\ V(s) \qquad \text{in } S^3;$$

$$y_\varepsilon(s) = \cosh\varepsilon\ x(s) + \sinh\varepsilon\ V(s) \qquad \text{in } H^3.$$

The vector field $V(s)$ then points "across" the ribbon.

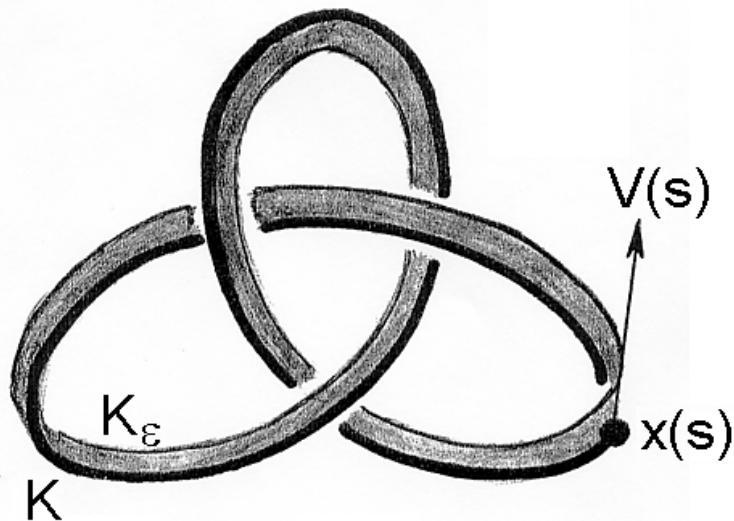

**THEOREM 5. LINK = TWIST + WRITHE in $S^3$ and $H^3$.**

**(1) On $S^3$ in left-translation format:**

$$\mathrm{Lk}(K, K_\varepsilon) = \mathrm{Tw}_L(V) + \mathrm{Wr}_L(K).$$

**(2) On $S^3$ in parallel transport format:**

$$\mathrm{Lk}(K, K_\varepsilon) = \mathrm{Tw}_P(V) + \mathrm{Wr}_P(K).$$

**(3) On $H^3$ in parallel transport format:**

$$\mathrm{Lk}(K, K_\varepsilon) = \mathrm{Tw}_P(V) + \mathrm{Wr}_P(K).$$

## 11. Proof scheme for Theorem 5.

In spirit, our proof of Theorem 5 for ribbons in $S^3$ and $H^3$ follows Calugareanu's original proof in $R^3$: we begin with the linking integrals given in Theorem 1 for the edges $K$ and $K_\varepsilon$ of our ribbon, let $\varepsilon$ shrink to zero, and observe the behavior of the linking integrand.

Of course, the value of the linking integral is independent of $\varepsilon$ for $\varepsilon > 0$, since the linking number is invariant under this kind of homotopy. But the linking integrand develops a $\delta$-function type singularity along the diagonal of $K \times K$, which is handled as follows.

Outside an appropriately chosen neighborhood of the diagonal, the linking integrand converges to the writhing integrand as $\varepsilon \to 0$, and its integral converges to the writhing number of the curve $K$. Inside this neighborhood of the diagonal, the linking integrand blows up, but its integral converges to the total twist of the normal vector field $V$ along $K$.

The crucial thing, recognized by Calugareanu, is that the width of the neighborhood of the diagonal in $K \times K$ must go to zero *much more slowly* than the width $\varepsilon$ of the ribbon. In fact, we will choose the neighborhood of the diagonal to have width $\varepsilon^p$, where $0 < p < 1/3$.

To give a sense of this in action, we will outline here the proof of Theorem 5, part (2), dealing with Link = Twist + Writhe in parallel transport format on $S^3$. The proofs for $H^3$ and for left-translation format on $S^3$ are essentially the same. In particular, in left-translation format, the integrand of the second integral in the expression for the linking number converges uniformly to the corresponding integrand for the writhing number.

Consider, in parallel transport format on $S^3$, the linking integrand of $K$ with $K_\varepsilon$,

$$F_\varepsilon(s, t) = (1/(4\pi^2))\, dx/ds \bullet P_{x(s) y_\varepsilon(t)}\big(dy_\varepsilon/dt \times \nabla_{y_\varepsilon(t)}\, \varphi(\alpha(x(s), y_\varepsilon(t)))\big)$$

and the writhing integrand of $K$,

$$F_0(s, t) = (1/(4\pi^2))\, dx/ds \bullet P_{x(s) x(t)}\big(dx/dt \times \nabla_{x(t)}\, \varphi(\alpha(x(s), x(t)))\big),$$

where $\varphi(\alpha) = (\pi - \alpha) \csc \alpha$.

Then the linking number of K and $K_\varepsilon$ is given by

$$Lk(K, K_\varepsilon) = \iint_{0 \leq s,t \leq L} F_\varepsilon(s, t) \, ds \, dt \, ,$$

and the writhing number of K is given by

$$Wr_P(K) = \iint_{0 \leq s,t \leq L} F_0(s, t) \, ds \, dt \, .$$

Since $\alpha(x, y)$ is the distance between x and y, and since $\varphi$ has a singularity just at $\alpha = 0$, the only difficulty in considering the convergence of the linking integral as $\varepsilon \to 0$ happens near the diagonal, where $s = t$.

To handle this, we first show that because the singularity of $\varphi''(\alpha)$ at $\alpha = 0$ is like $1/\alpha^3$, we have that for sufficiently small $\varepsilon > 0$,

$$|F_\varepsilon(s, t) - F_0(s, t)| \leq C \, \varepsilon^{1-3p} \, ,$$

provided that $|s - t| \geq \varepsilon^p$.

If $0 < p < 1/3$, then $1 - 3p > 0$, and hence $C \, \varepsilon^{1-3p} \to 0$ as $\varepsilon \to 0$. Therefore $F_\varepsilon(s, t)$ converges uniformly to $F_0(s, t)$ in the region $|s - t| \geq \varepsilon^p$. Since the writhing integrand $F_0(s, t)$ remains bounded even along the diagonal, this shows that

$$\iint_{|s-t| \geq \varepsilon^p} F_\varepsilon(s, t) \, ds \, dt \to \iint_{0 \leq s,t \leq L} F_0(s, t) \, ds \, dt = Wr_P(K) \, ,$$

as $\varepsilon \to 0$, that is, a *portion* of the linking integral converges to the *entire* writhing integral.

The more delicate part of the argument is the integral near the diagonal. A careful analysis using the Frenet equations reveals that for $0 < p < 1$,

$$\lim_{\varepsilon \to 0} \int_{t-\varepsilon^p}^{t+\varepsilon^p} F_\varepsilon(s, t) \, ds = (1/(2\pi)) \, x'(t) \times V(t) \bullet V'_P(t) \, .$$

Hence

$$\iint_{|s-t| \leq \varepsilon^p} F_\varepsilon(s, t) \, ds \, dt \to \int_{0 \leq t \leq L} (1/(2\pi)) \, x'(t) \times V(t) \bullet V'_P(t) \, dt = Tw_P(V) \, .$$

That is, the remaining portion of the linking integral converges to the entire twisting integral. In this way, we see that

$$Lk(K, K_\varepsilon) = Tw_P(V) + Wr_P(K) \, .$$

# 12. Helicity of vector fields on $S^3$ and $H^3$.

If $V$ is a smooth vector field on $R^3$ with compact support, then the Woltjer-Moffatt formula for its helicity,

$$H(V) = 1/(4\pi) \int_{R^3 \times R^3} V(x) \times V(y) \bullet (x-y)/|x-y|^3 \ dx \ dy \ ,$$

can be written succinctly as

$$H(V) = \int_{R^3} BS(V)(y) \bullet V(y) \ dy \ .$$

This is how Woltjer originally presented his invariant, $\int A \bullet B \ d(vol)$, with the role of $V$ played by the magnetic field $B$ and the role of $BS(V)$ played by its vector potential $A$. The helicity of a vector field, as we mentioned earlier, is a measure of the extent to which its orbits wrap and coil around one another.

We use the above formula to define the helicity of a vector field $V$ on $S^3$ or $H^3$, and then immediately obtain explicit integral formulas from Theorem 2.

**THEOREM 6. HELICITY INTEGRALS in $S^3$ and $H^3$.**

**(1)** *On $S^3$ in left-translation format:*

$$H(V) = \int_{S^3 \times S^3} L_{yx^{-1}} V(x) \times V(y) \bullet \nabla_y \varphi_0(x, y) \ dx \ dy$$

$$- 1/(4\pi^2) \int_{S^3 \times S^3} L_{yx^{-1}} V(x) \bullet V(y) \ dx \ dy$$

$$+ 2 \int_{S^3 \times S^3} \nabla_y \left( L_{yx^{-1}} V(x) \bullet \nabla_y \varphi_1(x, y) \right) \bullet V(y) \ dx \ dy \ ,$$

*where* $\varphi_0(\alpha) = (1/(4\pi^2))(\pi-\alpha) \cot \alpha$ *and* $\varphi_1(\alpha) = (-1/(16\pi^2)) \alpha (2\pi-\alpha)$.

**(2)** *On $S^3$ in parallel transport format:*

$$H(V) = \int_{S^3 \times S^3} P_{yx} V(x) \times V(y) \bullet \nabla_y \varphi(x, y) \ dx \ dy \ ,$$

*where* $\varphi(\alpha) = (1/(4\pi^2))(\pi - \alpha) \csc \alpha$.

**(3)** *On $H^3$ in parallel transport format:*

$$H(V) = \int_{H^3 \times H^3} P_{yx} V(x) \times V(y) \bullet \nabla_y \varphi(x, y) \ dx \ dy \ ,$$

*where* $\varphi(\alpha) = (1/(4\pi)) \operatorname{csch} \alpha$.

- In formula (1), if V is divergence-free, then the third integral in the definition of H(V) vanishes.

- In formulas (1) and (2), if the smooth vector field V on $S^3$ is divergence-free, then its helicity is the same as its *asymptotic* (or *mean*) *Hopf invariant*, as defined by Arnold (1974), and is invariant under the group of volume-preserving diffeomorphisms of $S^3$.

- In formula (3), we assume that V has compact support, in order to guarantee convergence of the integral.

- If V is a left-invariant unit vector field on $S^3$, recall that BS(V) = − ½ V. Then its helicity is given by

$$H(V) = \int_{S^3} BS(V)(y) \bullet V(y) \, dy$$

$$= \int_{S^3} -½ V(y) \bullet V(y) \, dy$$

$$= -½ \, vol(S^3) = -½ \, (2\pi^2)$$

$$= -\pi^2 .$$

If V is a right-invariant unit vector field on $S^3$, then BS(V) = ½ V, and hence

$$H(V) = +\pi^2 .$$

# 13. Upper bounds for helicity in $R^3$, $S^3$ and $H^3$.

Let $\Omega$ be a compact, smoothly bounded subdomain of $R^3$, $S^3$ or $H^3$, and V a smooth vector field defined on $\Omega$. Thinking of V as a current flow, its magnetic field BS(V) is defined by the same formulas as before, except that the integration is carried out only over $\Omega$.

If $\Omega \subset R^3$, then

$$BS(V)(y) = \int_\Omega V(x) \times \nabla_y \varphi(x, y) \, dx,$$

where $\varphi(\alpha) = -1/(4\pi\alpha)$.

If $\Omega \subset S^3$, then, for example in parallel transport format, we have

$$BS(V)(y) = \int_\Omega P_{yx}V(x) \times \nabla_y \varphi(x, y) \, dx,$$

where $\varphi(\alpha) = (-1/(4\pi^2))(\pi - \alpha) \csc \alpha$.

If $\Omega \subset H^3$, then we have

$$BS(V)(y) = \int_\Omega P_{yx}V(x) \times \nabla_y \varphi(x, y) \, dx,$$

where $\varphi(\alpha) = (-1/(4\pi)) \csch \alpha$.

The magnetic field BS(V) is defined throughout the ambient space. It is continuous everywhere, but its first derivatives suffer a discontinuity as one crosses the boundary of $\Omega$. This is a familiar situation from electrodynamics in Euclidean 3-space.

In what follows, we will restrict BS(V) to $\Omega$, and ignore its behavior outside this domain.

Let VF($\Omega$) denote the space of all smooth vector fields on $\Omega$, with the $L^2$ inner product

$$<V, W> = \int_\Omega V \bullet W \, d(\text{vol}),$$

and associated **energy** $<V, V>$ and **norm** $|V| = <V, V>^{1/2}$.

We seek a bound for the energy or norm of the output magnetic field BS(V) on $\Omega$ in terms of that of the input current flow V. Or to put it another way, we seek an upper bound for the $L^2$ operator norm of the Biot-Savart operator,

$$BS : VF(\Omega) \to VF(\Omega),$$

in terms of the geometry of the underlying domain $\Omega$.

As a consequence, we will determine an upper bound for the helicity $H(V) = <BS(V), V>$ of the vector field V in terms of its energy $<V, V>$ and the geometry of $\Omega$.

**THEOREM 7.** *Let $\Omega$ be a compact, smoothly bounded subdomain of $\mathbf{R}^3$, $\mathbf{S}^3$ or $\mathbf{H}^3$, and let $R = R(\Omega)$ be the radius of a ball in that space having the same volume as $\Omega$. Let V be a smooth vector field defined on $\Omega$. Then*

$$|\mathbf{BS(V)}| \leq \mathbf{N(R)} |\mathbf{V}|,$$

*where in $\mathbf{R}^3$ we have $\mathbf{N(R) = R}$;*

*in $\mathbf{S}^3$ we have $\mathbf{N(R) = (1/\pi) \big(2(1 - \cos R) + (\pi - R) \sin R\big)}$;*

*in $\mathbf{H}^3$ we have $\mathbf{N(R) = \sinh R}$.*

*It follows immediately that the helicity $\mathbf{H(V) = <BS(V), V>}$ is bounded by*

$$|\mathbf{H(V)}| \leq \mathbf{N(R)} |\mathbf{V}|^2.$$

The proof in $\mathbf{R}^3$ appears in our "upper bounds" paper (2000). Using the formulas for BS(V) provided in Theorem 2 above, the arguments in $\mathbf{S}^3$ and $\mathbf{H}^3$ are then substantially the same as in that paper.

Here are three observations.

• In $\mathbf{R}^3$, the overestimate $N(R) = R$ for the norm of BS grows like the *cube root* of the volume $(4/3)\pi R^3$ of $\Omega$.

• By contrast, in $\mathbf{H}^3$, the overestimate $N(R) = \sinh R$ for the norm of BS grows like the *square root* of the volume $2\pi(\sinh R \cosh R - R)$ of $\Omega$.

• In $\mathbf{S}^3$, the overestimate $N(\pi) = 4/\pi$ for the norm of BS on the entire 3-sphere can be compared with the known exact value, $|BS| = 1/2$.

# 14. Spectral geometry of the curl operator in $R^3$, $S^3$ and $H^3$.

Again let $\Omega$ be a compact, smoothly bounded subdomain of $R^3$, $S^3$ or $H^3$. We are interested here in curl eigenfields on $\Omega$, that is, vector fields $V$ on $\Omega$ which satisfy $\nabla \times V = \lambda V$, for $\lambda \neq 0$. In $R^3$ these eigenfields are used to model stable plasma flows; see our paper (1999).

Curl eigenfields exist for every value of $\lambda$. For example, in $R^3$, if

$$V = \sin \lambda z \; \mathbf{i} + \cos \lambda z \; \mathbf{j},$$

then
$$\nabla \times V = \lambda V.$$

We want to constrain the choice of vector fields $V$ by interior and boundary conditions which guarantee a discrete spectrum, while at the same time being reasonable for physical applications. We then seek a lower bound for the absolute values of the nonzero eigenvalues.

To begin, consider again the infinite-dimensional vector space $VF(\Omega)$ of smooth vector fields on $\Omega$, with the $L^2$ inner product.

Let $K(\Omega) \subset VF(\Omega)$ denote the subspace consisting of vector fields which are divergence-free and tangent to the boundary of $\Omega$,

$$K(\Omega) = \{V \in VF(\Omega) : \nabla \bullet V = 0, \; V \bullet \mathbf{n} = 0\},$$

where $\mathbf{n}$ denotes the unit outward normal vector field along the boundary $\partial\Omega$ of $\Omega$. We think of the vector fields in $K(\Omega)$ as *fluid knots* (hence the letter "K") because, in the traditional passage from geometric knot theory to fluid dynamics, a knot is modeled by an incompressible fluid flow, defined in a tubular neighborhood of the knot and tangent to its boundary.

To keep our problem manageable, we are only going to look for curl eigenfields within the subspace $K(\Omega)$ of $VF(\Omega)$.

Let $CK(\Omega) \subset K(\Omega)$ denote the subspace of vector fields whose curl lies in $K(\Omega)$. Thus

$$curl : CK(\Omega) \to K(\Omega).$$

For example, any divergence-free vector field on $\Omega$ which vanishes on $\partial\Omega$ lies in $K(\Omega)$. And clearly, any eigenfield of the curl operator in $K(\Omega)$ must actually lie in $CK(\Omega)$.

It is easy to see that a vector field $V \in K(\Omega)$ lies in the subspace $CK(\Omega)$ if and only if the circulation of $V$ around small loops on $\partial\Omega$ vanishes.

We want to know whether $curl : CK(\Omega) \to K(\Omega)$ is self-adjoint with respect to the $L^2$ inner product. In other words, is it true that

$$<\nabla \times V, W> = <V, \nabla \times W> ?$$

**Fact.** If $\Omega$ is simply connected, or equivalently, if all the components of $\partial\Omega$ are 2-spheres, then $curl : CK(\Omega) \to K(\Omega)$ is self-adjoint, but not otherwise.

So we seek further sensible boundary conditions which will make $curl$ self-adjoint, even when the domain $\Omega$ is not simply connected.

To this end, let $\Omega^*$ denote the closure of the complement of $\Omega$ in $R^3$, $S^3$ or $H^3$. Let g denote the total genus of $\partial\Omega$, that is, the sum of the genera of its components. Then, using real coefficients, $H_1(\partial\Omega)$ is a 2g-dimensional vector space, while $H_1(\Omega)$ and $H_1(\Omega^*)$ are each g-dimensional, and we have the direct sum

$$H_1(\partial\Omega) = \ker(H_1(\partial\Omega) \to H_1(\Omega)) \oplus \ker(H_1(\partial\Omega) \to H_1(\Omega^*)),$$

where the above homomorphisms are induced by the inclusions $\partial\Omega \subset \Omega$ and $\partial\Omega \subset \Omega^*$.

Let $AK(\Omega) \subset CK(\Omega) \subset K(\Omega)$ consist of all vector fields $V$ in $K(\Omega)$ whose circulation vanishes around any loop on $\partial\Omega$ which bounds in $\Omega^*$. We call $AK(\Omega)$ the space of **Amperian knots** because, by Ampere's Law, the magnetic field due to a current running entirely within $\Omega$ will have zero circulation around all loops on $\partial\Omega$ which bound in $\Omega^*$.

Magnetic fields generated by current distributions in $\Omega$ are divergence-free, but not in general tangent to the boundary of $\Omega$. So we define the ***modified Biot-Savart operator***

$$BS' : K(\Omega) \to K(\Omega)$$

by

$$BS'(V) = L^2 \text{ orthogonal projection of } BS(V) \text{ into } K(\Omega).$$

**Remark.** It is a simple consequence of the Hodge Decomposition Theorem that we have an $L^2$ orthogonal direct sum

$$VF(\Omega) = K(\Omega) \oplus G(\Omega),$$

where $G(\Omega)$ is the subspace of smooth gradient fields on $\Omega$. See, for example, our paper (2002). So, for $V$ in $K(\Omega)$,

$$BS(V) = BS'(V) + \text{ the gradient component of } BS(V).$$

**Some basic facts.**

- When $\Omega$ is simply connected, we have $AK(\Omega) = CK(\Omega)$, but when $\Omega$ is not simply connected, then $AK(\Omega)$ is smaller than $CK(\Omega)$.

- $AK(\Omega)$ is a dense subspace of $K(\Omega)$ in the $L^2$ topology.

- $AK(\Omega)$ is the image of $BS'$, and the map $BS' : K(\Omega) \to AK(\Omega)$ is a bijection.

- The maps

$$BS' : K(\Omega) \to AK(\Omega) \quad \text{and} \quad curl : AK(\Omega) \to K(\Omega)$$

are inverses of one another.

**THEOREM 8.** *Let $\Omega$ be a compact, smoothly bounded subdomain of $R^3$, $S^3$ or $H^3$, and let $R = R(\Omega)$ be the radius of a ball in that space having the same volume as $\Omega$.*

*Then* curl: $AK(\Omega) \to K(\Omega)$ *is a self-adjoint operator, and for each $V \in AK(\Omega)$ we have*

$$|\nabla \times V| \geq 1/N(R) \ |V|,$$

*where $N(R)$ is defined as in Theorem 7.*

*In particular, if $\lambda$ is the curl eigenvalue on $AK(\Omega)$ with smallest absolute value, then*

$$|\lambda| \geq 1/N(R).$$

The bounds given in Theorem 7 for the norm of BS(V) certainly hold also for its orthogonal projection BS'(V) into the subspace $K(\Omega)$.

Since *curl* : $AK(\Omega) \to K(\Omega)$ is the inverse of BS' : $K(\Omega) \to AK(\Omega)$, we simply invert these bounds to get Theorem 8.

University of Pennsylvania
Philadelphia, PA  19104

deturck@math.upenn.edu
gluck@math.upenn.edu